\documentclass[12pt,reqno]{amsart}
\usepackage{amsmath}%
\usepackage{amsfonts}%
\usepackage{amssymb}%
\usepackage[margin=.5in]{geometry}
\usepackage{subfig}
\usepackage{array,multirow}
\usepackage{ifpdf}
\ifpdf
\usepackage{graphicx}
\else
\usepackage[draft]{graphicx}
\fi
\usepackage{color}
\usepackage{enumerate}
\usepackage[T1]{fontenc}
\usepackage{caption}
\usepackage{booktabs,multirow}
\usepackage{enumerate}
\usepackage{listings}
\usepackage{color}
\usepackage{caption}
\usepackage{xcolor}
\usepackage{hyperref}
\hypersetup{%
  colorlinks=false,
  linkbordercolor=red,
  pdfborderstyle={/S/U/W 1}
}

\definecolor{dkgreen}{rgb}{0,0.6,0}
\definecolor{gray}{rgb}{0.5,0.5,0.5}
\definecolor{mauve}{rgb}{0.58,0,0.82}

\newcommand{\bint}{\displaystyle{\int\hspace{-11.4pt}\mathit{8}}}
\AtBeginDocument{}

\begin{document}

\title[]{On the branching geometry of algebraic functions}
\author{Dominic C. Milioto}
\email{icorone@hotmail.com}
\date{\today}
\begin{abstract}
This paper describes an algorithm for determining the branching geometry of algebraic functions.  The graphs of these complex-valued functions have a complicated  interweaving structure that can be described by analytic branches separated by singular points.  Power expansions for the branches in discs centered at a point can be computed using the Newton Polygon method, and expansions around annular regions centered at the origin computed using a version of Laurent's Theorem applied to algebraic functions.  However, neither of the methods enable a determination of the region of convergence of the power series.  In this paper, a method using analytic continuation is used to determine the domain of analyticity for the branches, and the Root Test used to numerically check the results.  
 
The software used to implement the algorithm is Mathematica ver. 11.1.
\end{abstract}
\maketitle
\section{Introduction}
The objects studied in this paper are algebraic functions $w(z)$ expressed implicitly by the equation
\begin{equation}
f(z,w)=a_0(z)+a_1(z)w+a_2(z)w^2+\cdots+a_{n}(z)w^{n}=0,
\label{eqn001}
\end{equation}
with $z$ and $w$ complex variables and the coefficients, $a_i(z)$, polynomials in $z$ with rational coefficients.  The degree of the function is the highest power of $w$.   
 By the Implicit Function Theorem, this equation defines locally, an analytic function $w(z)$ when $\displaystyle \frac{\partial f}{\partial w}\neq 0$.  And it is known from the general theory of algebraic functions that the solution, $w(z)$, has $n$ fractional power series expansions about a point. These fractional power series can be computed by the Newton polygon method and often have finite radii of convergences determined by the nearest impinging singular point so do not in general represent the entire function.  Power series for annular regions are also fractional and can be computed by a version of Laurent's Theorem applied to algebraic functions.
 
 The purpose of this paper is four-fold:
\begin{enumerate}
\item{Describe an algorithm to determine the branching geometry of an algebraic function},
\item{Describe a method of applying Laurent's Theorem to compute annular power series expansions of an algebraic function}, 
\item{Determine the region of convergence of a fractional power series representation of $w(z)$,} 
\item {Present a software tool for visually investigating the branching of algebraic functions.}
\end{enumerate}
\section{Some properties of algebraic functions used in this paper}
Fractional power expansions of algebraic functions are called Puiseux series and usually have finite radii of convergences.  As stated earlier, series centered at a point can be computed by the Newton Polygon method but usually represent a small portion of (\ref{eqn001}).  Power series for annular regions can in principle, be computed using a variation of Laurent's Theorem for multivalued functions.  The Laurent expansions computed this way are again fractional power expansions.  In this way, the global branching geometry of $w(z)$ can be represented in the form of singular points segregating the z-plane into annular regions where the function is analytic and ramifies into multivalued branches.  This branching, both around the singular points and annular regions, are represented by the notation $\{s_1,s_2,\cdots,s_n\}$ where each $s_i$ represents an analytic and single-valued sheet of the branch.  For example the notation $\{1,2\}$ represents a $2$-cycle branch, such as $\sqrt{z}$,  and the numbers 1 and 2 represent the sort order of the function value at predetermined points in the z-plane called the annular or singular point reference point.  Figure \ref{figure:plot1} shows these points as the blue and red dots labeled `arp' and `srp'.

 The Resultant of $f(z,w)$ with it's derivative $f_w$ is denoted by $R(f,f_w)$.  A point $z$ where $R(f,f_w)$ is zero is a singular point of $f$.  However this does not tell us which branch sheet is singular.  A point where $a_n(z)=0$ is a pole, possibly ramified.  Singular points are sorted first by the real component and then the imaginary component and labeled $s_n$ with $n$ ranging from one to the total number of singular points including the point at infinity.  Even though the function may not ramify and thus not be singular at infinity we still include it in the list of singular points.

 A singular point may not affect all coverings of an algebraic function unless the function fully-ramifies at the singular point.  For example, a $10$-degree function may only ramify into a single $2$-cycle branch at a singular point with the remaining coverings single-cycled and unramified.  In this case, the $2$-cycle covering is singular.  The eight single-cycle coverings are not analytically affected at this singular point unless one is affected by a pole of the function.  However, if the function fully-ramifies into a $10$-cycle branch at this singular point, all coverings would be affected.  It is for this reason a power expansion of an algebraic function often has a region of convergence extending beyond the nearest singular point:  the branch coverings may simply not be affected by the singular point.  Only when the covering becomes singular does the convergence region of its power expansion become established. The main objective of this paper is to identify which singular point is interrupting the analyticity of branch cycles thereby establishing the region of convergence of their power expansions.
\begin{figure}
	\centering
			\includegraphics[scale=.5]{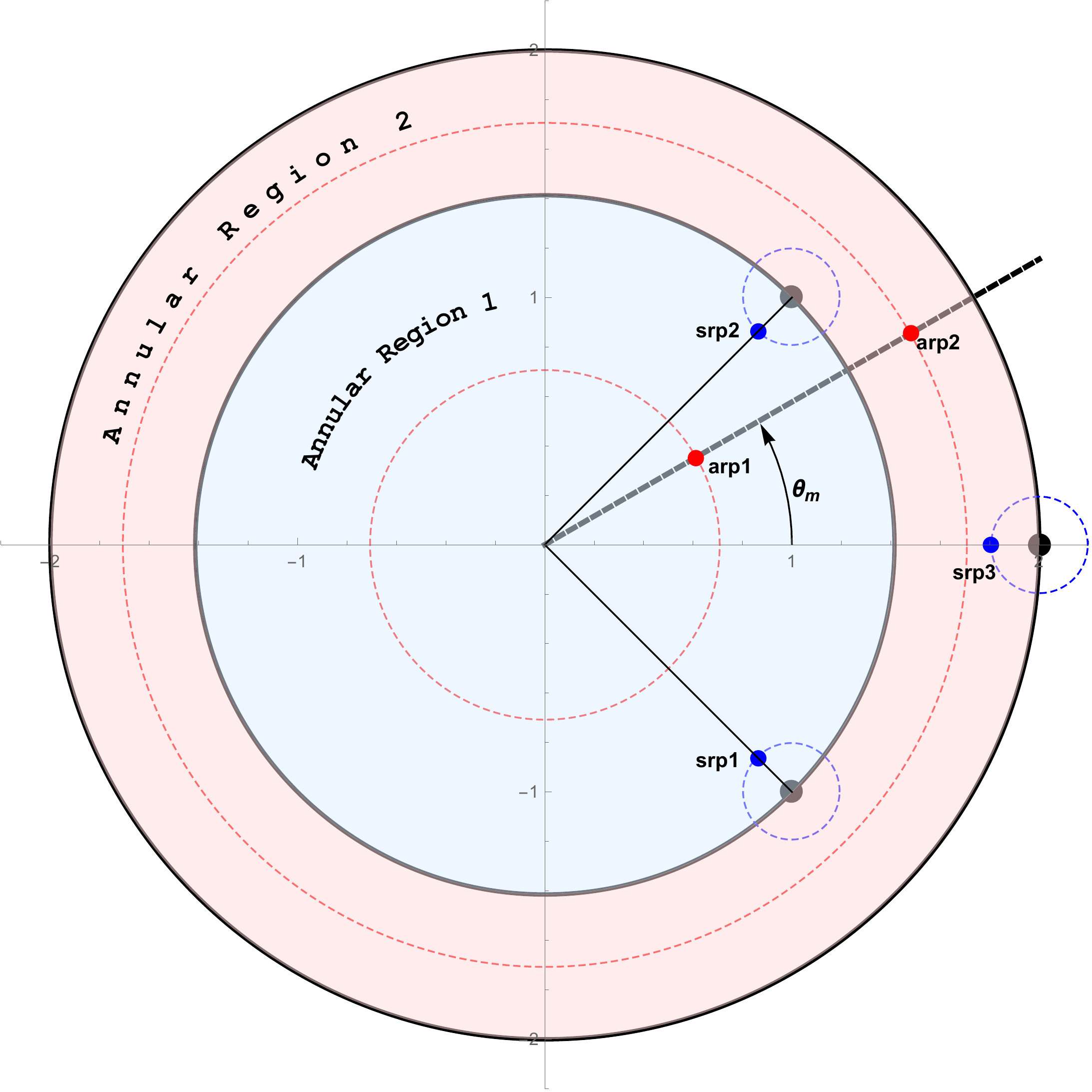}
		\caption{Annular Regions and Integration Paths}
   \label{figure:plot1}
\end{figure}

 The following conventions are used in this paper:
\begin{enumerate}
\item The Puiseux expansions of $(\ref{eqn001})$ consist of a set of $d$-valued branches.  A branch is sometimes called a $d$-cycle where $d$ is a positive integer.  A power series in $z^{1/3}$ would be a $3$-cycle branch.  It has three coverings over the complex $z$-plane similar to $\sqrt[3]{z}$.  The sum of the cycles is always equal to the degree of the function in $w$.
\item The concept of ``branch'' is used throughout this paper and refers to a multi-valued $d$-cycle of $w(z)$. 
\item The following discussion makes use of the term, ``extending a branch over a singular point''.  This is in reference to the discussion above about singularities and the coverings they affect.
\item $r_c$ is a positive integer representing the radius of a circular ring with center at the origin.  These rings are created by arranging the finite singular points in order of increasing absolute value.  The smallest non-zero singular point is therefore on ring one with radius $r_1$.  Singular points with the next largest absolute value are on ring two of radius of $r_2$ and so forth.  
\item We further sub-divide singular points into function singular points and branch singular points.  The function singular points is the solution set to $R(f,f_w)=0$ as discussed above.  A branch singular point is a singular point impinging one or more branch sheets, and this is done to emphasize singular points often do not affect all branches.  For example, if singular point $s_i$ is on a 3-cycle branch, then $s_i$ becomes the branch singular point.  Other cycles may not be affected by this singular point and so $s_i$ is not a branch singular point for these cycles.
\item Power series for algebraic branches are fractional power series called Puiseux series.  Two methods are used to compute these power series:  Newton Polygon and Laurent's Theorem applied to algebraic functions.  For brevity, only a few terms of a series are listed.  Actual computations for this paper used up to $1024$ terms with over $4000$ digits of precision.
\end{enumerate}
\section{Determining Monodromies}
 We first compute the finite singular points of the function by setting the resultant, $R(f,f_w)$, equal to zero.  This gives a list of $z$ values where the singular points are located.  We then sort the points in order of absolute value creating rings around the origin.  In the annular regions between rings, the function is analytic and splits into multi-valued branches. 

 Our next objective is to determine the branching geometry around each singular point and each annular region.  We then determine the annular domain of analyticity for each annular branch.  The branching is called monodromy or `ramified covering' and we use these terms interchangeably.   The monodromy is computed by numerical integration:  If we follow an analytically-continuous $2\pi$ circular path over a branch covering around a singular point enclosing no other singular point, then at the end of the path, we may or may not return to the same starting point on the branch.  If we return to the same value, we have identified a single-cycle branch.  If we need to traverse a $4\pi$ circuit to return to the same starting point, then we have traversed a $2$-cycle branch and so forth.  In practice however, we use a somewhat different algorithm but the principle is similar.  This give us the local monodromy around the singular point.  And we use the same approach to determine the monodromy of each annular region.    Naturally, this requires a comparison of floating-point numbers and some type of numerical accuracy must be established to identify terms considered identical.  We check a certain number of digits to the right of the decimal place without rounding.  For example, we could require five digits of accuracy between $w(0)$ and $w(2\pi)$ to be considered equal.  However care must be take to avoid rounding or cases when we have for example $0.999999999$.  We can check the later case and convert it to 1 in order to check the values. 

 Consider first the function
\begin{equation}
f(z,w)=w^2-(z-(1+i))^2(z-(1-i))^2(z-2)=0.
\label{equation:equation5}
\end{equation}
Setting $R(f,f_w)=0$ we obtain as the singular points, $\{1-i,1+i,2,\infty\}$ and the rings in Figure \ref{figure:plot1}. The black circles demarcate the annular regions separated by one or more singular points shown in the figure as black dots.  During the analysis, we integrate around annular regions and singular points beginning at reference points labeled 'arp' in the figure for annular regions and 'srp' for singular points.  The annular integration paths are shown as the red dashed circles.  The blue dashed circles are the singular point integration paths. 

 As stated, we use numerical integration over $2\pi$ routes around each annular region to determine the monodromy.  And in order to minimize errors we take a route with a mean radius of each region to maximize the distance to the nearest singular point.  For example, if a region was given by $1\leq r\leq 2$, we would integrate at 1.5 using $z(t)=1.5 e^{it}$.  But how do we decide if values are identical?  We can do this by extracting a set number of decimal digits.  For example, consider the value $1.23901950708 - 0.278993512 i$.  Taking an accuracy of 5 digits, we construct the following approximate integer representation of this number:
 $$
((1, (1, 2, 3, 9,0), 1), (-1, (2, 7, 8, 9,9), 0)),
$$
with the format $\{s, \{d_1, d_2, \cdots, d_n\}, e\}$ with $s$ being the sign of the number, $(d_1,d_2,\cdots,d_n)$ the digits of the number and $e$ being the number of digits to the left of the decimal place. Notice that we did not round the number.  So we basically convert a floating point number to an integer sequence so that we can make exact comparisons of floating point numbers up to a desired accuracy.  In practice, we are comparing a number against a set of numbers and we check that the smallest difference between the set of numbers is larger than the accuracy of the comparison to minimize errors.

 Now consider a 10-degree algebraic function.  We wish to integrate over one of the annular regions and then compare the starting and ending values of the function over each $2\pi$ route of each branch sheet.  Table \ref{table:table1} gives actual results for a 10-degree function.  Looking carefully at the starting and ending values we see the $w_1$ root goes back to the $w_1$ root so this is a single cycle branch $\{1\}$.  Next, the  $w_2$ root goes to $w_7$ then $w_6$ then $w_4$ back to $w_2$.  This is a $4$-cycle branch, $\{2,7,6,4\}$.  We then have two single cycle branches $w_3$ and $w_5$, then a $3$-cycle branch $w_8$, $w_9$, $w_{10}$.  Therefore the annular monodromy is
$$\{\{1\},\{2,7,6,4\},\{3\},\{5\},\{8,9,10\}\}.$$
\begin{table}[ht]
\caption{Monodromy digits} 
\centering 
\begin{tabular}{|c|c|c|} 
\hline\hline 
Root & Starting Value & Ending Value \\[0.5ex]
\hline
1 &(((-1,(1,0,5,0,8),1),((-1,(4,0,8,8,9),-1)) &((-1,(1,0,5,0,8),1),(-1,(4,0,8,8,9),-1)) \\
2 & (((-1,(4,7,6,5,7),0),(-1,(1,9,8,2,8),-1))& ((1,(5,8,2,9,7),0),(1,(9,7,6,1,3),0)) \\
3 &(((-1,(4,0,9,0,9),0),(1,(7,3,5,4,2),0)) &((-1,(4,0,9,0,9),0),(1,(7,3,5,4,2),0))  \\
4 &(((-1,(2,5,8,4,2),0),(-1,(1,3,8,2,8),1)) & ((-1,(4,7,6,5,7),0),(-1,(1,9,8,2,8),-1))\\
5 &(((-1,(2,3,9,4,1),0),(-1,(7,6,4,2,7),0)) &((-1,(2,3,9,4,1),0),(-1,(7,6,4,2,7),0))\\
6 &(((-1,(3,8,0,2,3),-1),(1,(7,0,9,5,3),-1)) &((-1,(2,5,8,4,2),0),(-1,(1,3,8,2,8),1))\\
7 &(((1,(5,8,2,9,7),0),(1,(9,7,6,1,3),0)) &((-1,(3,8,0,2,3),-1),(1,(7,0,9,5,3),-1))\\
8 &(((1,(9,5,0,3,7),0),(-1,(4,4,9,0,4),0)) &((1,(1,1,0,7,5),1),(-1,(5,2,0,1,7),-1)) \\
9 &(((1,(1,1,0,7,5),1),(-1,(5,2,0,1,7),-1)) &((1,(1,1,4,4,9),1),(1,(6,2,4,2,6),0)) \\
10 &(((1,(1,1,4,4,9),1),(1,(6,2,4,2,6),0)) &((1,(9,5,0,3,7),0),(-1,(4,4,9,0,4),0)) \\
\hline 
\end{tabular}
\label{table:table1} 
\end{table}
\section{Numerically Solving the Monodromy Differential Equations}
 In order to determine monodromies, we integrate over the function around each annular region and around each singular point.  We do this by solving the monodromy differential equation:
Given $f(z,w)=a_0(z)+a_1(z)w+\cdots+a_n(z)w^n=0$, we have $\displaystyle \frac{dw}{dt}=-\frac{f_z}{f_w}\frac{dz}{dt}$ which is a first-order differential equation.  In order to compute the monodromy, we can solve this differential equation for $w(t)$ over a circular path around annular regions and singular points.  For example, we would solve the following $n$ initial value problems for a function of degree $n$:
$$
\frac{dw}{dt}=-\frac{f_z}{f_w}\frac{dz}{dt}, \; (z_0,w_i),\; i=1,2,\cdots,n
$$
where $w_i$ are the roots to the expression $f(z_0,w)=0$.  We then compile the values of $w(t)$ at the beginning and ending of each $2\pi$ route. 
 However, in order to obtain accurate results, we will have to adjust the working precision and step size of the numerical integration as needed to achieve the desired accuracy.  To economize this, we set up a loop of decreasing step size with increasing precision ranging from $(1/1000,20)$ to $(1/50000,65)$.  However, even at the high range of precision, we may not be able to resolve the branching if for example the region is very small or the branch sheets are very close to each other.  For random functions, we have found if the annular size, $|r_2-r_1|$,  is larger than $1/5000$, and the difference between branch sheets is greater than the desired accuracy, then a 10-degree function can usually be processed successfully.  Therefore, this paper deals only with functions of degree 10 or less with $|r_2-r_1| \geq 1/5000$.
\section{Determining the branching geometry}
We divide this section into the following sub-sections:
\begin{enumerate}
\item{Computation of singular points,}
\item{Construction of annular regions,}
\item{Compute singular point monodromies,}
\item{Compute annular monodromies,}
\item{Determine branch-continuation support,}
\item{Compute branch continuations.}
\item{Determine continuations over poles},
\end{enumerate}
and use the following function to illustrate the concepts:
\begin{equation}
f(z,w)=(9)+(-7)w+(7-2 z-4 z^2-z^3)w^2+(7 z-2 z^3)w^3.
\label{equation:equation10}
\end{equation}   
\subsection{Computation of singular points}
\begin{table}[ht]
\caption{Singular points} 
\centering 
\begin{tabular}{|c|c|} 
\hline 
Singularity & Value\\[0.5ex]
\hline
 $s_1$ & \color{red}0. \\
 $s_2$ & \color{blue}-0.21713-0.255535 i  \\
 $s_3$ & \color{blue}-0.21713+0.255535 i \\
 $s_4$ & \color{blue}1.36571 \\
 $s_5$ & \color{blue}-1.83037-0.0204249 i  \\
 $s_6$ & \color{blue}-1.83037+0.0204249 i \\
 $s_7$ & \color{red}-1.87083 \\
 $s_8$ & \color{red}1.87083 \\
 $s_9$ & \color{blue}2.21932\, -1.22168 i \\
 $s_{10}$ & \color{blue}2.21932\, +1.22168 i \\
 $s_{11}$ & \color{blue}-6.85468-6.05671 i  \\
 $s_{12}$ & \color{blue}-6.85468+6.05671 i \\
 $s_{13}$ & $\infty$  \\
\hline 
\end{tabular}
\label{table:table3} 
\end{table}
 Using the built-in Mathematica function \texttt{NSolve}, the finite singular points are easily determined.  However, in order to later compute successfully the annuli and singular monodromies, we need to compute the singular points at a sufficiently high precision.  In this paper, we do so at 75 digits of precision and we attempt to carry approximately, this precision throughout the calculations. Table \ref{table:table3} lists the singular points in blue where we have identified poles in red. 

\begin{table}[ht]
\caption{Rings and Singular points} 
\centering 
\begin{tabular}{|c|c|c|c|>{$}c<{$}|}\hline
    \text{Ring} & Radius & \text{Singularity} & \text{Value} \\\hline
    \multirow{2}{*}{$r_1$} & \multirow{2}{*}{0.335327} &  $s_1$ & \color{blue}-0.21713 - 0.255535 i \\
                           &                           &  $s_2$ & \color{blue}-0.21713 + 0.255535 i \\ \hline
                     $r_2$ &                   1.36571 &  $s_3$ &              \color{blue}1.36571 \\ \hline
    \multirow{2}{*}{$r_3$} &  \multirow{2}{*}{1.83048} &  $s_4$ & \color{blue}-1.83037 - 0.0204249 i \\
                           &                           & $s_5$ &  \color{blue}-1.83037 + 0.0204249 i \\ \hline
    \multirow{2}{*}{$r_4$} &  \multirow{2}{*}{1.87083} & $s_6$ &           \color{red}-1.87083 \\
                           &                           & $s_7$ &               \color{red}1.87083 \\ \hline
    \multirow{2}{*}{$r_5$} &  \multirow{2}{*}{2.53336} &    $s_8$ & \color{blue}2.21932 - 1.22168 i \\
                           &                           &    $s_9$ & \color{blue}2.21932 + 1.22168 i \\ \hline
    \multirow{2}{*}{$r_6$} &  \multirow{2}{*}{9.14715} &   $s_{10}$ & \color{blue}-6.85468 - 6.05671 i \\
                           &                           &   $s_{11}$ & \color{blue}-6.85468 + 6.05671 i \\ \hline
  \end{tabular}
	\label{table:table59} 
\end{table}
\subsection{Construction of rings and annular regions}
 Once the singular points are computed, they are arranged in order of increasing absolute values which then segregate the z-plane into annular regions separated by rings.  On each ring lies one or more singular points.  Table \ref{table:table59} lists the rings, $r_i,$ which will become the regions of convergence of the power expansions of each function branch computed below.    For example, a branch may have a region of convergence given by $(r_1,r_3)$.

 Between each ring are annular region devoid of singular points. Table \ref{table:table4} lists the annular regions separated by the ring singular points.  The last region is simply determined by an arbitrary distance from the most distant finite singular point.  In this case, this distance is $4$.
\subsection{Computation of annular monodromies}

\begin{table}[ht]
\caption{Annular Regions} 
\centering 
\begin{tabular}{|c|c|} 
\hline
Annulus & Annulus/singular point \\[0.5ex]
\hline 
   & \color{red}0.  \\
 $c_1$ & \{0.00005,0.335277\} \\
   & \color{blue}-0.21713-0.255535 i  \\
   & \color{blue}-0.21713+0.255535 i \\
 $c_2$ & \{0.335377,1.36566\}  \\
   & \color{blue}1.36571  \\
 $c_3$ & \{1.36576,1.83043\}  \\
   & \color{blue}-1.83037-0.0204249 i \\
   & \color{blue}-1.83037+0.0204249 i \\
 $c_4$ & \{1.83053,1.87078\} \\
   & \color{red}-1.87083 \\
   & \color{red}1.87083  \\
 $c_5$ & \{1.87088,2.53331\}  \\
   & \color{blue}2.21932\, -1.22168 i \\
   & \color{blue}2.21932\, +1.22168 i \\
 $c_6$ & \{2.53341,9.1471\}  \\
   & \color{blue}-6.85468-6.05671 i  \\
   & \color{blue}-6.85468+6.05671 i \\
 $c_7$ & \{9.1472,13.1471\}  \\
   & $\infty$  \\ 
 \hline 
\end{tabular}
\label{table:table4} 
\end{table}

 With annular regions defined, we next compute the annular monodromies.  One way to compute this is to simply integrate around annuli over $2\pi$ routes, one route for each branch sheet for a total of $n$ sheets and determine how many routes to return to a starting point as was shown above in Table \ref{table:table1} until all the branch sheets have been processed.  This necessarily involves comparing floating-point numbers, but if we are willing to accept a tolerance say of five decimal digits to the right of the decimal point or other numerical accuracy, experience has shown we can obtain reliable results.
 As stated earlier, to effect this integration, we integrate the monodromy differential equation.   In the case of the annular regions, we choose a circular path midway in the region with starting value $z(\theta_0)$ with $\theta_0$ chosen somewhat arbitrarily to best effect the integration.  For example, $\theta$ could be chosen to maximize the distance from most singular points.  We then form a table of starting and ending values as described in Section 3 above for each path and then determine the monodromies through integer comparisons.  This gives us the monodromy results in Table \ref{table:table5}.
\begin{table}[ht]
\caption{Annular Monodromies} 
\centering 
\begin{tabular}{|c|c|c|} 
\hline
Annulus & Annulus/singular point &Monodromies\\[0.5ex]
\hline 
   & \color{red}0.  & \\
 $c_1$ & \{0.00005,0.335277\} & ((1),(2),(3)) \\
   & \color{blue}-0.21713-0.255535 i &\\
   & \color{blue}-0.21713+0.255535 i &\\
 $c_2$ & \{0.335377,1.36566\} & ((1,3,2)) \\
   & \color{blue}1.36571 &  \\
 $c_3$ & \{1.36576,1.83043\} & ((1,3),(2)) \\
   & \color{blue}-1.83037-0.0204249 i &\\
   & \color{blue}-1.83037+0.0204249 i &\\
 $c_4$ & \{1.83053,1.87078\} & ((1,2),(3))\\
   & \color{red}-1.87083 &\\
   & \color{red}1.87083 & \\
 $c_5$ & \{1.87088,2.53331\} & ((1,2),(3)) \\
   & \color{blue}2.21932\, -1.22168 i &\\
   & \color{blue}2.21932\, +1.22168 i &\\
 $c_6$ & \{2.53341,9.1471\}  & ((1,2),(3))\\
   & \color{blue}-6.85468-6.05671 i  &\\
   & \color{blue}-6.85468+6.05671 i &\\
 $c_7$ & \{9.1472,13.1471\} & ((1),(2,3)) \\
   & $\infty$ &\\ 
 \hline 
\end{tabular}
\label{table:table5} 
\end{table}
\subsection{Computing singular point monodromies}
The procedure for computing singular point monodromies is identical to that for the annular regions:  We solve the monodromy DE for a circular path around each singular point containing no other singular points.  In the case of the monodromy around infinity, we integrate over a closed contour enclosing all finite singular points.  These results are shown in Table \ref{table:table6}. 
\begin{table}[ht]
\caption{Annular and Singular Point Monodromies} 
\centering 
\begin{tabular}{|c|c|c|} 
\hline
Annulus & Annulus/singular point &Monodromies\\[0.5ex]
\hline 
   & \color{red}0. & \color{red}((1),(2),(3)) \\
 $c_1$ & \{0.00005,0.335277\} & ((1),(2),(3)) \\
   & \color{blue}-0.21713-0.255535 i & \color{blue}((1),(2,3)) \\
   & \color{blue}-0.21713+0.255535 i & \color{blue}((1),(2,3))\\
 $c_2$ & \{0.335377,1.36566\} & ((1,3,2)) \\
   & \color{blue}1.36571  & \color{blue}((1),(2,3))\\
 $c_3$ & \{1.36576,1.83043\} & ((1,3),(2)) \\
   & \color{blue}-1.83037-0.0204249 i & \color{blue}((1),(2,3))\\
   & \color{blue}-1.83037+0.0204249 i &\color{blue}((1),(2,3))\\
 $c_4$ & \{1.83053,1.87078\} & ((1,2),(3))\\
   & \color{red}-1.87083 & \color{red}((1),(2),(3)\\
   & \color{red}1.87083 & \color{red}((1),(2),(3)) \\
 $c_5$ & \{1.87088,2.53331\} & ((1,2),(3)) \\
   & \color{blue}2.21932\, -1.22168 i & \color{blue}((1,2),(3))\\
   & \color{blue}2.21932\, +1.22168 i & \color{blue}((1,2),(3))\\
 $c_6$ & \{2.53341,9.1471\}  & ((1,2),(3))\\
   & \color{blue}-6.85468-6.05671 i  & \color{blue}((1,2),(3))\\
   & \color{blue}-6.85468+6.05671 i & \color{blue}((1,2),(3))\\
 $c_7$ & \{9.1472,13.1471\} & ((1),(2,3)) \\
   & $\infty$ & \color{blue}((1),(2,3)) \\ 
 \hline 
\end{tabular}
\label{table:table6} 
\end{table}
\subsection{Determining support}
 Now that we have the monodromies for both the annular regions and singular points, we can determine the possible support of each annular branch across intervening singular points.  A necessary condition for branch continuation across a singular point is that that singular monodromy between successive annuli must support a sufficient number of single-cycle branches to continue the branch into the next annular region.  For example, if we are considering extending a 3-cycle branch across a singular point, then the monodromy around this singular point must have three single-cycle branches.  And likewise for other cycles.

 A second necessary condition is that the next or post-annular monodromy must have the same branch cycle type we are attempting to continue.   Table \ref{table:table7} shows the possible branch continuations.
\begin{table}[ht]
\caption{Support} 
\centering 
\begin{tabular}{|c|c|c|c|} 
\hline
Annulus & Annulus/singular point &Monodromies & Support\\[0.5ex]
\hline 
   & \color{red}0. & \color{red}((1),(2),(3))& \\
 $c_1$ & \{0.00005,0.335277\} & ((1),(2),(3)) &\\
   & \color{blue}-0.21713-0.255535 i & \color{blue}((1),(2,3)) &\\
   & \color{blue}-0.21713+0.255535 i & \color{blue}((1),(2,3))&\\
 $c_2$ & \{0.335377,1.36566\} & ((1,3,2)) &\\
   & \color{blue}1.36571  & \color{blue}((1),(2,3)) & \\
 $c_3$ & \{1.36576,1.83043\} & ((1,3),(2))  & ((2))  \\
   & \color{blue}-1.83037-0.0204249 i & \color{blue}((1),(2,3)) & \\
   & \color{blue}-1.83037+0.0204249 i &\color{blue}((1),(2,3)) & \\
 $c_4$ & \{1.83053,1.87078\} & ((1,2),(3)) & ((1,2),(3))  \\
   & \color{red}-1.87083 & \color{red}((1),(2),(3) & \\
   & \color{red}1.87083 & \color{red}((1),(2),(3)) & \\
 $c_5$ & \{1.87088,2.53331\} & ((1,2),(3)) & ((3))  \\
   & \color{blue}2.21932\, -1.22168 i & \color{blue}((1,2),(3)) & \\
   & \color{blue}2.21932\, +1.22168 i & \color{blue}((1,2),(3)) & \\
 $c_6$ & \{2.53341,9.1471\}  & ((1,2),(3)) & ((3))  \\
   & \color{blue}-6.85468-6.05671 i  & \color{blue}((1,2),(3)) & \\
   & \color{blue}-6.85468+6.05671 i & \color{blue}((1,2),(3)) & \\
 $c_7$ & \{9.1472,13.1471\} & ((1),(2,3)) & \\
   & $\infty$ & \color{blue}((1),(2,3))  & \\ 
 \hline 
\end{tabular}
\label{table:table7} 
\end{table}
For example, consider the first and second annuli in Table \ref{table:table7} and the intervening singular points:  The first annulus has three single-cycle branches.  In order to continue one or more of these branches into annulus two, the singular points between these regions must support single-cycle branches.  In this case, the singular points do have single-cycle branches.  However, annulus 2 is fully-ramified into a 3-cycle branch so does not support continuing any of the single-cycle branches.  And likewise for annulus 2:  In order to continue this branch into annulus 3, we would need three single-cycle branches over the intervening singular points.  In the case of annulus 3 with monodromy $\{\{1,3\},\{2\}\}$, we cannot continue the $2$-cycle branch into annulus 4 since the intervening singular points do not support two single-cycle branches.  However there is support to continue branch $\{2\}$:  the singular point has a single-cycle branch and annulus 4 does as well.  Now consider annuli 4 and 5 and the two poles between them:  the poles do not ramify but rather consists of three single-cycle branches and so supports continuing both the $\{1,2\}$ branch and the $\{3\}$ branch however not holomorphically:  one or more of the branches will contain a pole.  
\subsection{Determining branch continuations}
In the previous section, we defined a necessary condition for branch continuation:  The intervening singular points must support continuation.  A sufficient condition is that each branch sheet be analytically continuous over all singular points on the bordering ring and continuing into the next annular region onto a branch with the same monodromy sequence.  For example, if we are attempting to continue a 3-cycle branch $\{2,1,3\}$  across a singular point from annulus $k$ to annulus $k+1$, then from above, the singular monodromy must support three single-cycle branches, and annulus $k+1$ must have a $3$-cycle branch such as $\{1,2,3\}$ such that sheet $2$ in annulus $k$ is continued onto sheet $1$ of annulus $k+1$, sheet $1$ is continued onto sheet $2$, and sheet $3$ is continued onto sheet $3$ or: 
$$
\begin{array}{cccc} 2 & 1 & 3 \\
                   \downarrow & \downarrow & \downarrow \\
									 1 & 2 & 3
									\end{array}
$$

 After we have determined the possible continuations, we can use numerical integration to determine if each branch covering is analytically continuous around the singular points. 
In order to show how this integration is effected, we use Equation (\ref{equation:equation5}) again and the plot shown in Figure \ref{figure:figure8}.
\begin{figure}
\centering
\includegraphics[scale=0.5]{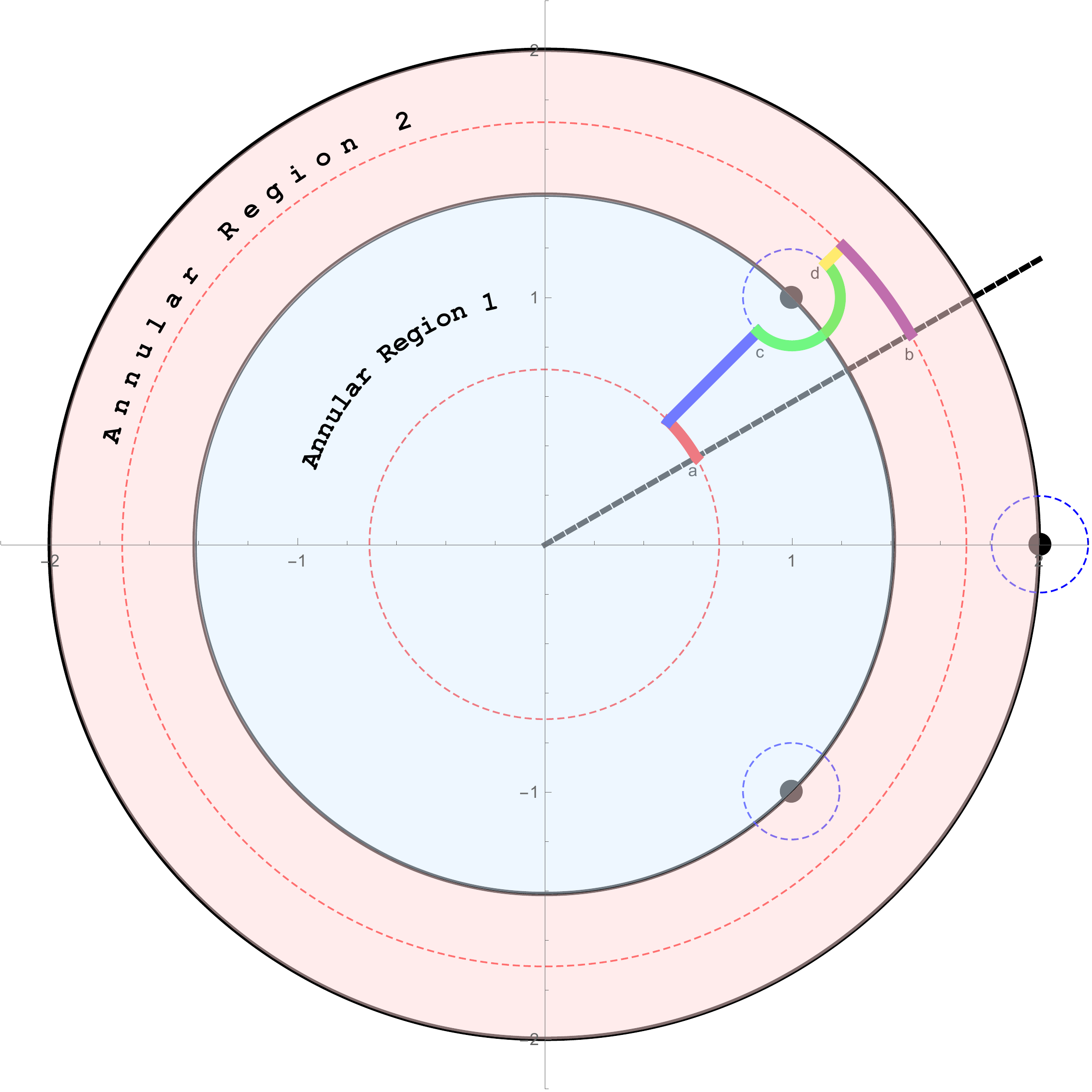}
\caption{Continuation Path}
\label{figure:figure8}
\end{figure}
In order to check if a branch sheet in region 1 is continued onto a branch sheet in region 2 across the singular point at $1+i$, we analytically continue the branch sheet to region 2 by integrating from point a to b in Figure \ref{figure:figure8} along the red, blue, green, yellow and purple contours.  The dotted black line is the argument where all annular monodromies are taken at their associated reference point, arp.  At point a in the diagram, the annular monodromy of the annulus 1 is determined.  At point $b$, the monodromy of  annulus $2$ is determined.  In this way we can compare the  monodromies between annulus 1 and 2 to determine if a branch in region 1 analytically continues to a branch in region 2.  However we cannot simply integrate from point a to b to determine this as we may be integrating too close to another singular point. Rather, we integrate over the third green leg, around the singular point to point $d$ where we choose the radius of the green leg such that it minimizes the distance to the nearest singular point while remaining in the bordering rings.  We now have integrated over an analytically-continuous route from region 1 to region 2 to point $d$.  And  our objective is to get to point $b$ where we can compare annular monodromies.  And although we could in principle integrate directly from point $d$ to point $b$, we take a more symmetrical route first along the short yellow leg and then the purple leg to $b$.
At the point $b$ in the figure, we can compare the monodromies of region 2 to those in region 1.  The resulting continuations are shown in the fourth column of Figure \ref{table:continuationtable1} with arrows between branch regions signifying continuations.
\begin{figure}
\centering
\includegraphics[scale=0.6]{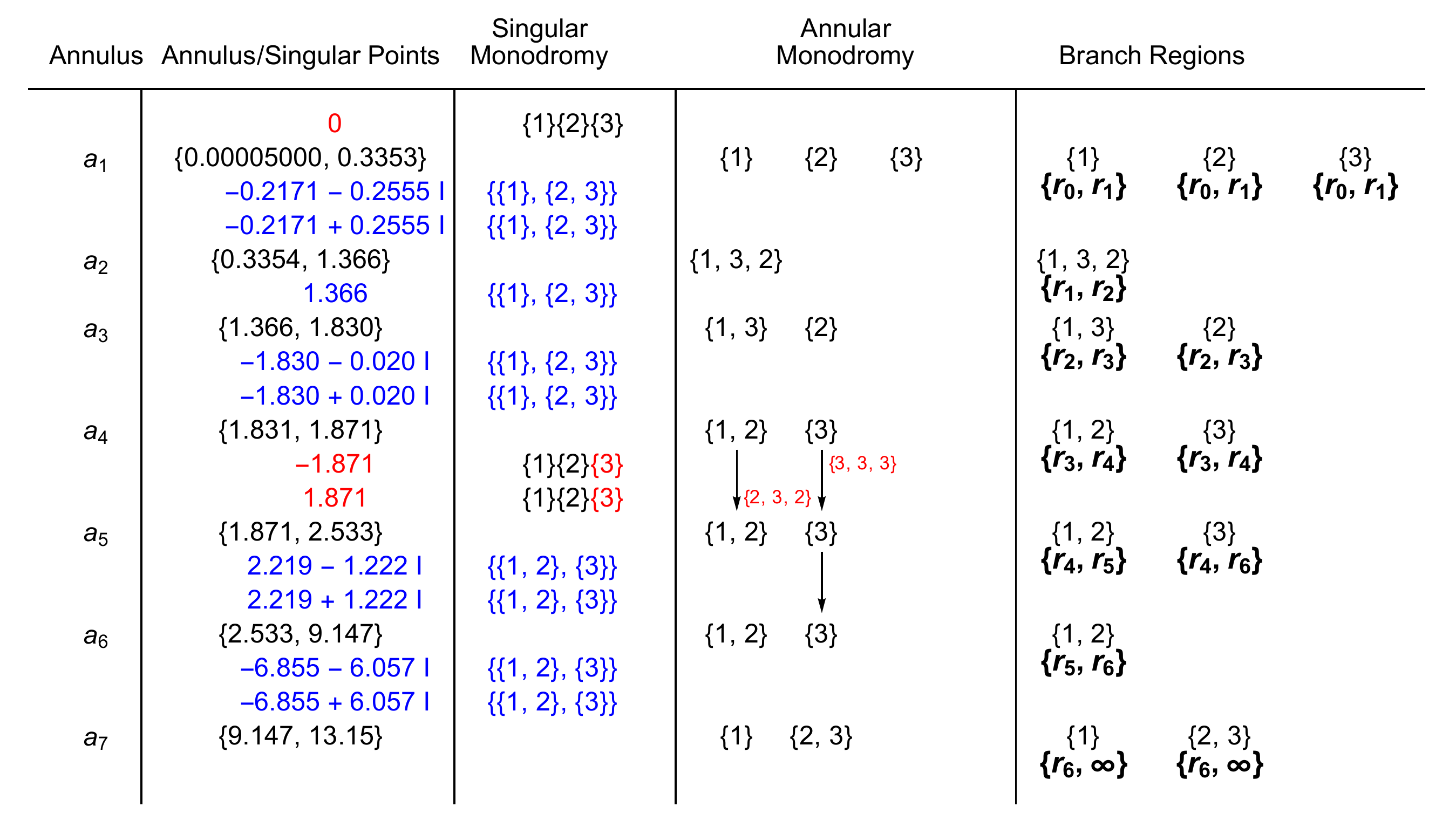}
\caption{Branch Regions for Equation \ref{equation:equation10}}
\label{table:continuationtable1}
\end{figure}   
\subsection{Continuations over poles}
The arrows between branches in annuli 4 and 5 in Figure \ref{table:continuationtable1} shows the $\{1,2\}$ and $\{3\}$ cycles in ring 4 continued into ring 5 across two singular points.  However since the intervening singular points are poles, one or more of the branch continuations may be meromorphic, that is, contain a pole. 
 In cases involving poles, we can determine which branch sheet is affected by the pole by computing the Puiseux expansion around each singular point and then computing the value of each series at the singularity reference point $c$ in Figure \ref{figure:figure8}.  We first compute the the function at point $c$ and sort the values, then continue each branch sheet to point $c$ and determine for each sheet, the sort order at point $c$.  We then compute the Puiseux series at point $c$ and then compute the value of each series at point $c$.  Note however, the point at $c$ for the Puiseux expansions will be offset by the singular point.  For each singular point, we record the annulus, branch, sheet, singular point, sort index, the value of $z$ at the singular point reference point, and it's offset, $zP$, when we translate axes to the singular point.  These results are shown in Table \ref{table:table10}.
\begin{table}[ht]
\caption{Singular Witness Data} 
\centering 
\begin{tabular}{|c|c|c|c|c|c|c|c|} 
\hline
 \text{Line Num} &\text{Annulus} & \text{Next Cycle} & \text{Sheet} & \text{sing pt} & \text{Sort index} & \text{zstart} & \text{zP} \\
\hline\hline 
1 & 4 & \{1,2\} & 1 & -1.87083 &  2  & -1.85741 & 0.0134149 \\
2 & 4 & \{1,2\} & 2 & -1.87083 & 1  & -1.85741 & 0.0134149 \\
3 & 4 & \{1,2\} & 1 & 1.87083 & 1 & 1.85741 & -0.0134149 \\
4 & 4 & \{1,2\} & 2 & 1.87083 & 3  & 1.85741 & -0.0134149 \\
5 & 4 & \{3\} & 3 & -1.87083 & 3 & -1.85741 & 0.0134149 \\
6 & 4 & \{3\} & 3 & 1.87083 & 2 & 1.85741 & -0.0134149 \\
\hline
\end{tabular}
\label{table:table10} 
\end{table}
Now look at the first line in the table:  We are checking the first sheet of the 2-cycle branch $\{1,2\}$ in annulus 4 across singular point $-1.87083$ and we find when we continue this branch to the singular point reference point $c$ in Diagram \ref{figure:figure8} , we find it continues onto the second root of the function at this point.  The value of $z$ at point $c$ is $-1.85741$ and it's offset value around the singular point is $0.0134149$.  We next compute the Puiseux expansions at this singular point:
$$
\begin{aligned}
w_1(z)&=(1.06397-1.26645 i)+(3.01595-6.43854 i) z+(0.4263-73.5684 i) z^2\\
&-(426.834+907.601 i) z^3-(12030.5+9669.1 i) z^4-(249021.+55181. i) z^5+\cdots\\
w_2(z)&=(1.06397+1.26645 i)+(3.01595+6.43854 i) z+(0.4263+73.5684 i) z^2\\
&-(426.834-907.601 i) z^3-(12030.5-9669.1 i) z^4-(249021.-55181. i) z^5+\cdots\\
w_3(z)&=-1.76336+\frac{0.234968}{z}-5.65797 z-0.676369 z^2+853.755 z^3+24061.1 z^4+\cdots\\
\end{aligned}
$$
and then compute the value of each series at the offset point $zP=0.0134149$:
$$
\begin{aligned}
w_1(zP)&=1.10295\, -1.36858 i \\
w_2(zP)&=1.10295\, +1.36858 i \\
w_3(zP)&= 15.6793 \\
\end{aligned}
$$
Now, when we compute the sorted list of the function at the point $c$ we obtain:
$$
\begin{aligned}
1.10295 -1.36858 i\\
1.10295 +1.36858 i\\
15.6793\\
\end{aligned}
$$
and note $w_3(z)$ is the pole and it's value at point $c$ is approx $15.6793$ and that is the third index into the sorted list of function values at point $c$.  Since we determined that the first sheet of branch $\{1,2\}$ was continued onto the second index, then we know this sheet is not affected by the pole.  And in the second line of the table, we see the second sheet of this branch is continued onto index 1.  Therefore, since branch $\{1,2\}$ is not continued onto the pole sheet, we know this 2-cycle branch is not affected by this singular point.

Now look at lines 3 and 4 in Table \ref{table:table10} as we attempt to continue the $\{1,2\}$ cycle across $1.87083$.  We find at point $c$, the branches continue onto indexes 1 and 3.  At the singular reference point for this pole (the equivalent of point $c$ for the pole), we again compute the Puiseux series at $1.87083$:
$$
\begin{aligned}
w_1(z)&=-0.951784+0.494400 z-0.146001 z^2+0.0672361 z^3-0.0564447 z^4+0.0251513 z^5+\cdots \\
w_2(z)&=0.546915-0.405452 z+0.375163 z^2-0.435910 z^3+0.584467 z^4-0.853222 z^5+\cdots\\
w_3(z)&=-0.566856-\frac{1.23497}{z}+0.179982 z-0.377397 z^2+0.449107 z^3-0.571336 z^4+\cdots
\end{aligned}
$$
and note $w_3(z)$ is the pole.  When we compute the values of the series at the singular reference point we obtain:
$$
\begin{aligned}
w_1(zP)&=-0.958443\\
w_2(zP)&=0.552423 \\
w_3(zP)&=91.4904.
\end{aligned}
$$
And the sorted list of function values at the witness mark is:
$$
\begin{aligned}
-0.958443\\
0.552423\\
91.4904\\
\end{aligned},
$$
so that index 3 is the pole.  Thus we find the second sheet of this branch continues across this singular point but in a meromorphic fashion.  When we analyze the single cycle branch in the same manner, we find it continues meromorphically across the negative pole.  We identify these continuations across singular points in Table \ref{table:table10} as the small red monodromies along the continuation arrows. In the first case we have the second sheet of $(1,2)$ was continuous onto the third sheet of the negative pole to the second sheet of the post annulus.  This is the $(2,3,2)$ on the side of the first arrow.  And since the single cycle branch $\{3\}$ was meromorphically continued onto the positive pole to $\{3\}$ on the post annulus, we label this as $(3,3,3)$ next to the continuation arrow between these branches.  So that all branches in this annulus are affected by one of the poles.  Therefore the annular region of convergence of the branches in $a_4$ is contained in the annulus.  This is shown as $(r_4,r_5)$ below the branch monodromies in the table.  And the only branch which extends beyond an annular region is branch $\{3\}$ in $a_5$ where it is shown to have a region of convergence $(r_4,r_6)$.  To give a concrete example of this, the annular Puiseux series for this branch is
$$
\{3\}(z)=\cdots+\frac{1.8811}{z^3}-\frac{0.9475}{z^2}+\frac{1.3012}{z}+0.0072 z-0.00053 z^2+0.000031 z^3+\cdots
$$
and has a region of convergence of $(r_5,r_7)$ or approximately $(1.87088,9.1471)$.  

\section{Computation of annular Puiseux series}
The Newton Polygon algorithm computes power series for the function centered at a point.  These power series have radii of convergences equal to the absolute value of the branch singular point which often extends only into the first few rings.  In order to compute power series for the function in the remaining annular rings, we use a version of Laurent's Theorem applied to algebraic functions: 
\begin{equation}
\begin{aligned}
w_n(z)&= \sum_{k=0}^{\infty} a_k (z^{1/n})^k+\sum_{k=1}^{\infty} \frac{b_k}{\left(z^{1/n}\right)^k} \\
&=A(z)+S(z) \\
a_k&=\frac{1}{2n\pi i} \bint \frac{w_n(z)}{\left(z^{1/n}\right)^{k+n}} dz\\
b_k&=\frac{1}{2n\pi i} \bint w_n(z)\left(z^{1/n}\right)^{k-n} dz.\\
\label{equation:equation4}
\end{aligned}
\end{equation}
with $A(z)$ being the analytic terms of the series and $S(z)$, the singualr terms.  
Or in symmetrical form:
\begin{equation}
\begin{aligned}
a_k&=\frac{1}{2n\pi i} \bint \frac{w_n(z)}{\left(z^{1/n}\right)^{k+n}} dz\\
w_n(z)&=\sum_{p=-\infty}^\infty a_p (z^{1/n})^p,
\end{aligned}
\end{equation}
where the integral symbol $\bint$ indicates the integration is over a closed analytically continuous route along the integrand branch surfaces.   For example, if the integrand contained a 4-cycle branch, the integration would proceed over the branch surface along an analytically continuous $8\pi$ route of winding number 4.  However, we cannot simply integrate the expression as we would encounter implied branch-cuts.  We can avoid this by the following derivation:

\begin{equation}
\begin{aligned}
c_k&=\frac{1}{2n\pi i} \bint \frac{w_n(t)}{\left(z^{1/n}\right)^{k+n}} dz=\frac{1}{2 n \pi} \int_{t_0}^{t_e} \frac{w_n(t)re^{it}}{\left(re^{it}\right)^{\frac{k+n}{n}}}\\
&=\frac{1}{2 n\pi} \int_{t_0}^{t_e} w_n(t) \left(re^{it}\right)^{-k/n} dt\\
&=\frac{1}{2 n\pi r^{k/n}} \int_{t_0}^{t_e} w_n(t)\big[\cos(tk/n)-i\sin(tk/n)\big] dt\\
&=\frac{1}{2 n\pi r^{k/n}} I(k,n).
\label{eqn5}
\end{aligned}
\end{equation}

From the Root Test the upper radius of convergence of the analytic terms is
$$R_a=\frac{\displaystyle 1}{\displaystyle\limsup_{k\to \infty}\; a_k^{k/n}}$$
and the lower radius of convergence of the singular terms is
$$R_s=\limsup_{k\to -\infty} \;a_k^{k/n}.$$  And in order to estimate these limits numerically, we plot $\biggr(1/k,r\biggr(\frac{2 n \pi}{|I(k,n)|}\biggr)^{n/k}\biggr),\; k \neq 0$ and extrapolate to $1/k \to 0$.  This will necessarily invert the expression for the upper limit so that we write the region of convergence as

\begin{equation}
R=
\biggr\{\displaystyle\limsup_{k\to -\infty}\;r\biggr(\frac{2 n \pi}{|I(k,n)|}\biggr)^{n/k} \; , \;
\displaystyle\liminf_{k\to\infty}\;r\biggr(\frac{2 n \pi}{|I(k,n)|}\biggr)^{n/k}\biggr\}.
\label{equation1000}
\end{equation}

\begin{figure}
	\centering
			\includegraphics[scale=.6]{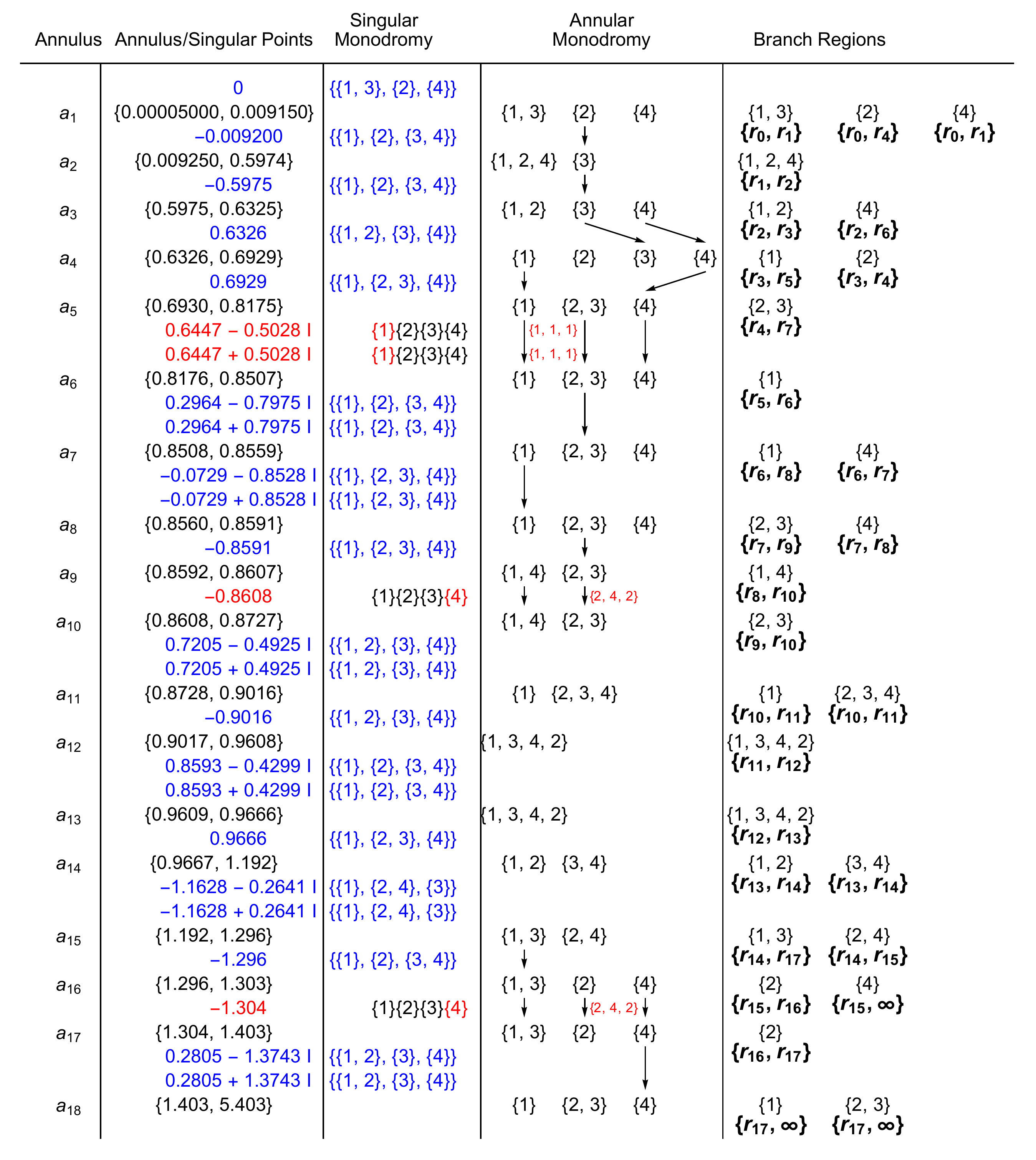}
		\caption{Branch Table for Equation \ref{equation:equation20}}
   \label{figure:continuationtable2}
\end{figure}

In order to demonstrate these formulas, we use
\begin{equation}
\begin{aligned}
  f(z,w)&=(-z^2+z^3)\\
	&+(-4 z+3 z^2)w\\
	&+(-z^3-9 z^4)w\\
	&+(-2+8 z+4 z^2-4 z^3)w^3\\
	&+(6-8 z^2+7 z^3+8 z^4)w^4.
	\end{aligned}
	\label{equation:equation20}
	\end{equation}
	The ring and singular points are shown in Table \ref{table:table52} and the continuations are in Figure \ref{figure:continuationtable2}.
	
	\begin{table}[ht]
\caption{Rings and Singular points} 
\centering 
\begin{tabular}{|c|c|c|c|>{$}c<{$}|}
\hline
    \text{Ring} & Radius & \text{Singularity} & \text{Value} \\
		\hline\hline
$r_{1}$ & 0.00919971 & $s_{1}$ & \color{blue}-0.00919971\\ \hline
$r_{2}$ & 0.597463 & $s_{2}$ & \color{blue}-0.597463\\ \hline
$r_{3}$ & 0.632598 & $s_{3}$ & \color{blue}0.632598\\ \hline
$r_{4}$ & 0.692915 & $s_{4}$ & \color{blue}0.692915\\ \hline
\multirow{2}{*}{$r_{5}$} &\multirow{2}{*}{0.81757} & $s_{5}$ & \color{red}0.644655 - 0.502832 i\\
& & $s_{6}$ & \color{red}0.644655 + 0.502832 i\\ \hline
\multirow{2}{*}{$r_{6}$} &\multirow{2}{*}{0.85077} & $s_{7}$ & \color{blue}0.296412 - 0.797464 i\\
& & $s_{8}$ & \color{blue}0.296412 + 0.797464 i\\ \hline
\multirow{2}{*}{$r_{7}$} &\multirow{2}{*}{0.855943} & $s_{9}$ & \color{blue}-0.0728759 - 0.852835 i\\
& & $s_{10}$ & \color{blue}-0.0728759 + 0.852835 i\\  \hline

$r_{8}$ & 0.859144 & $s_{11}$ & \color{blue}-0.859144\\ \hline
$r_{9}$ & 0.86077 & $s_{12}$ & \color{red}-0.86077\\ \hline

\multirow{2}{*}{$r_{10}$} &\multirow{2}{*}{0.87273} & $s_{13}$ & \color{blue}0.72046 - 0.492539 i\\
& & $s_{14}$ & \color{blue}0.72046 + 0.492539 i\\ \hline

$r_{11}$ & 0.901619 & $s_{15}$ & \color{blue}-0.901619\\ \hline

\multirow{2}{*}{$r_{12}$} &\multirow{2}{*}{0.960847} & $s_{16}$ & \color{blue}0.859329 - 0.429862 i\\
& & $s_{17}$ & \color{blue}0.859329 + 0.429862 i\\ \hline

$r_{13}$ & 0.966603 & $s_{18}$ & \color{blue}0.966603\\ \hline

\multirow{2}{*}{$r_{14}$} &\multirow{2}{*}{1.19237} & $s_{19}$ & \color{blue}-1.16276 - 0.264081 i\\
& & $s_{20}$ & \color{blue}-1.16276 + 0.264081 i\\ \hline

$r_{15}$ & 1.29612 & $s_{21}$ & \color{blue}-1.29612\\ \hline
$r_{16}$  & 1.30354 & $s_{22}$ & \color{red}-1.30354\\ \hline
\multirow{2}{*}{$r_{17}$} &\multirow{2}{*}{1.4026} & $s_{23}$ & \color{blue}0.280488 - 1.37427 i\\
& & $s_{24}$ & \color{blue}0.280488 + 1.37427 i\\

\hline
 \end{tabular}
	\label{table:table52} 
\end{table}

	Consider now the first annular region in Figure \ref{figure:continuationtable2} having a 2-cycle and two single cycle branches.  Using the Newton-Polygon algorithm, we compute (numerically), the power series representation for these branches:
	\begin{equation}
	\begin{aligned}
	\{1,3\}(z)&=-1.4142 i z^{1/2}-2.8750 z+13.3300 i z^{3/2}+83.9648 z^2-611.2916 i z^{5/2}+\cdots \\
\{2\}(z)&= -0.2500 z+0.0703 z^2+0.02075 z^3+0.00399 z^4-0.1061 z^5+\cdots\\
\{4\}(z)&=0.3333+4.6667 z-168.22 z^2+9523.5 z^3-65812 z^4+5.056*10^7 z^5.
\end{aligned}
\end{equation}
From the branch table, we know the $\{1,3\}$ cycle in $a_1$ has a radius of convergences of $r_1$ while the $\{2\}$ branch extends to $r_4$. These are labeled below the branches in the Branch columns. 
   
 Once we have the continuations, we can update Table \ref{figure:continuationtable2} with the branch regions in column five.  These lists the branch and the annular region of convergence of the associated power expansion.  For example, in the monodromy column, we see $\{2\}$ continues from the first annulus to the fourth.  Therefore, the region column does not list the intervening continuations of this branch in the second third and fourth annular regions as these are part of the same $\{2\}$ branch.   In the region column for this branch we have the notation $\{r_0,r_4\}$.  This gives the annular region of convergence of the Puiseux expansion in terms of the ring numbers $\{r_1,r_5\}$.  And from Figure \ref{figure:continuationtable2}, we see the annular domain of convergence is approximately $\{0,0.692915\}$ which are computed from the absolute value of singular points. Likewise, the $\{4\}$branch in annulus 3 has a region of convergence of $(r_2,r_6)$ as it continues to the sixth annulus.  This function therefore has 29 branches where we do not treat a meromorphic continuation across a pole as a single branch but rather two distinct branches. 

\section{Using the Root Test to confirm the results}

In order to confirm independently the convergence data in Table \ref{figure:continuationtable2}, Equation \ref{eqn5} was used to compute power expansions for each annular branch.  Recall the symmetrical expression for the coefficients:
$$ c_n=\frac{I(k,n)}{2 n \pi r^{k/n}}$$
or

\begin{equation}
\frac{c_n}{2\pi n}=\frac{I(k,n)}{r^{k/n}}
\label{equation:eqn5002}.
\end{equation}

And consider the $(1,2,4)$ branch in annulus $2$ with a domain $(0.009,0.59)$ and the computation for the $a_{100}$ term using two different values for the radius of integration and the subsequent powers:  

$$
\begin{aligned}
0.5^{100/3}&=9\times 10^{-11}\\
0.01^{100/3}&=2\times 10^{-67}.
\end{aligned}
$$

Now, the left side of Equation \ref{equation:eqn5002} is a constant.  So in order for the expressions to hold for $r=0.01$ and $r=0.5$, the integral $I(100,3)$ would have to be much smaller for $r=0.01$ than when $r=0.5$, in fact, on the order of $10^{-50}$.  This would require a much more precise evaluation of the integral of the analytic terms when $r=0.01$.  So that we use a value of $r$ close to the upper limit of $0.5$ for the analytic terms in this case.  A similar argument holds for the computation of singular terms where we have
\begin{equation}
\frac{c_n}{2\pi n}=I(k,n)r^{|k|/n}
\label{equation:eqn50020}
\end{equation}
so that we use a value of $r$ close to the lower limit of $0.009$. 

 Additionally, when the annular region extends beyond the unit circle, we are then confronted with the opposite case.  Consider again the two terms:
\begin{equation}
\begin{aligned}
1.4^{100}=4\times 10^{14} \\
1.6^{100}=2\times 10^{20}.
\end{aligned}
\end{equation}

So only a small difference results in a difference of $10^{6}$.  So that when the domains are larger than the unit circle, we use a value of $r$ close to the lower limit for both the analytic and singular terms. 

And for regions smaller than the unit cirle and when the upper and lower boundaries of the domain are similar in magnitude, calculations for the analytic power terms used a radius of integration  close to the upper convergence boundary, and calculations for singular terms, a radius close to the lower boundary.  In practice, these rules were adjusted by trial-and-error for best results.

We wish to numerically estimate the integral $I(k,n)$ as accurately as reasonably possible.  However in many cases, the numeric estimate includes a small residual.  It is therefore important to visually inspect the results and confirm the difference between the actual value of the integral and the residual is as large as possible to minimize errors.  

\begin{table}[ht]
\caption{Root Test Results for Singular Series} 
\centering 
\begin{tabular}{|c| c| c| c|c|c| c|c|} 
\hline
 \text{Annulus} &\text{Branch} & \text{Radius} & \text{Precision} & \text{Range} & \text{Actual } & \text{Est.}& \text{$\%$ Error} \\
\hline\hline 

 2 & \{1,2,4\} & 0.0100 & 30 & \{-100,-200\} & 0.009250 & 0.009162 & 0.95 \\
 3 & \{1,2\} & 0.6100 & 30 & \{-50,-100\} & 0.5975 & 0.5931 & 0.74 \\
 3 & \{4\} & 0.6100 & 30 & \{-50,-100\} & 0.5975 & 0.5953 & 0.36 \\
 4 & \{1\} & 0.6511 & 30 & \{-100,-200\} & 0.6326 & 0.6291 & 0.56 \\
 4 & \{2\} & 0.6387 & 30 & \{-50,-100\} & 0.6326 & 0.6303 & 0.36 \\
 5 & \{2,3\} & 0.7093 & 30 & \{-50,-100\} & 0.6930 & 0.6882 & 0.69 \\
 6 & \{1\} & 0.8209 & 30 & \{-50,-100\} & 0.8176 & 0.8174 & 0.024 \\
 7 & \{1\} & 0.8516 & 30 & \{-200,-400\} & 0.8508 & 0.8500 & 0.091 \\
 7 & \{4\} & 0.8513 & 30 & \{-200,-300\} & 0.8508 & 0.8452 & 0.66 \\
 8 & \{2,3\} & 0.8565 & 30 & \{-200,-400\} & 0.8560 & 0.8585 & 0.29 \\
 8 & \{4\} & 0.8563 & 30 & \{-200,-300\} & 0.8560 & 0.8502 & 0.67 \\
 9 & \{1,4\} & 0.8605 & 50 & \{-300,-600\} & 0.8592 & 0.8480 & 1.3 \\
 10 & \{2,3\} & 0.8620 & 30 & \{-200,-400\} & 0.8608 & 0.8527 & 0.95 \\
 11 & \{1\} & 0.8757 & 30 & \{-200,-300\} & 0.8728 & 0.8670 & 0.66 \\
 11 & \{2,3,4\} & 0.8757 & 50 & \{-300,-600\} & 0.8728 & 0.8608 & 1.4 \\
 12 & \{1,3,4,2\} & 0.9076 & 30 & \{-200,-400\} & 0.9017 & 0.8938 & 0.87 \\
 13 & \{1,3,4,2\} & 0.9615 & 50 & \{-300,-600\} & 0.9609 & 0.9468 & 1.5 \\
 14 & \{1,2\} & 0.9892 & 30 & \{-200,-400\} & 0.9667 & 0.9754 & 0.91 \\
 14 & \{3,4\} & 0.9892 & 55 & \{-400,-600\} & 0.9667 & 0.959 & 0.81 \\
 15 & \{1,3\} & 1.195 & 50 & \{-300,-600\} & 1.192 & 1.183 & 0.75 \\
 15 & \{2,4\} & 1.193 & 50 & \{-300,-600\} & 1.192 & 1.183 & 0.75 \\
 16 & \{2\} & 1.296 & 30 & \{-50,-100\} & 1.296 & 1.293 & 0.28 \\
 16 & \{4\} & 1.337 & 30 & \{-100,-200\} & 1.296 & 1.294 & 0.17 \\
 17 & \{2\} & 1.305 & 30 & \{-50,-100\} & 1.304 & 1.300 & 0.29 \\
 18 & \{1\} & 1.443 & 30 & \{-200,-300\} & 1.403 & 1.394 & 0.61 \\
 18 & \{2,3\} & 1.443 & 50 & \{-300,-600\} & 1.403 & 1.392 & 0.78 \\
\hline
\end{tabular}
\label{table:singRootTest} 
\end{table}
 Generally between $50$ and $400$ terms of each series was computed with a precision between $30$ and $50$ digits.  The Root Test as per Equation \ref{equation1000} was then used to estimate the convergence domains by extrapolating the data to the point $1/n\to 0$ using the Mathematica build-in function \texttt{FindFormula}. Both Tables \ref{table:singRootTest} and \ref{table:analRootTest} show the percent error of the calculations as well as parameters used in the computations.

\begin{figure}[!ht]
     \subfloat[\label{ring2branches:1}]{%
       \includegraphics[width=0.4\textwidth]{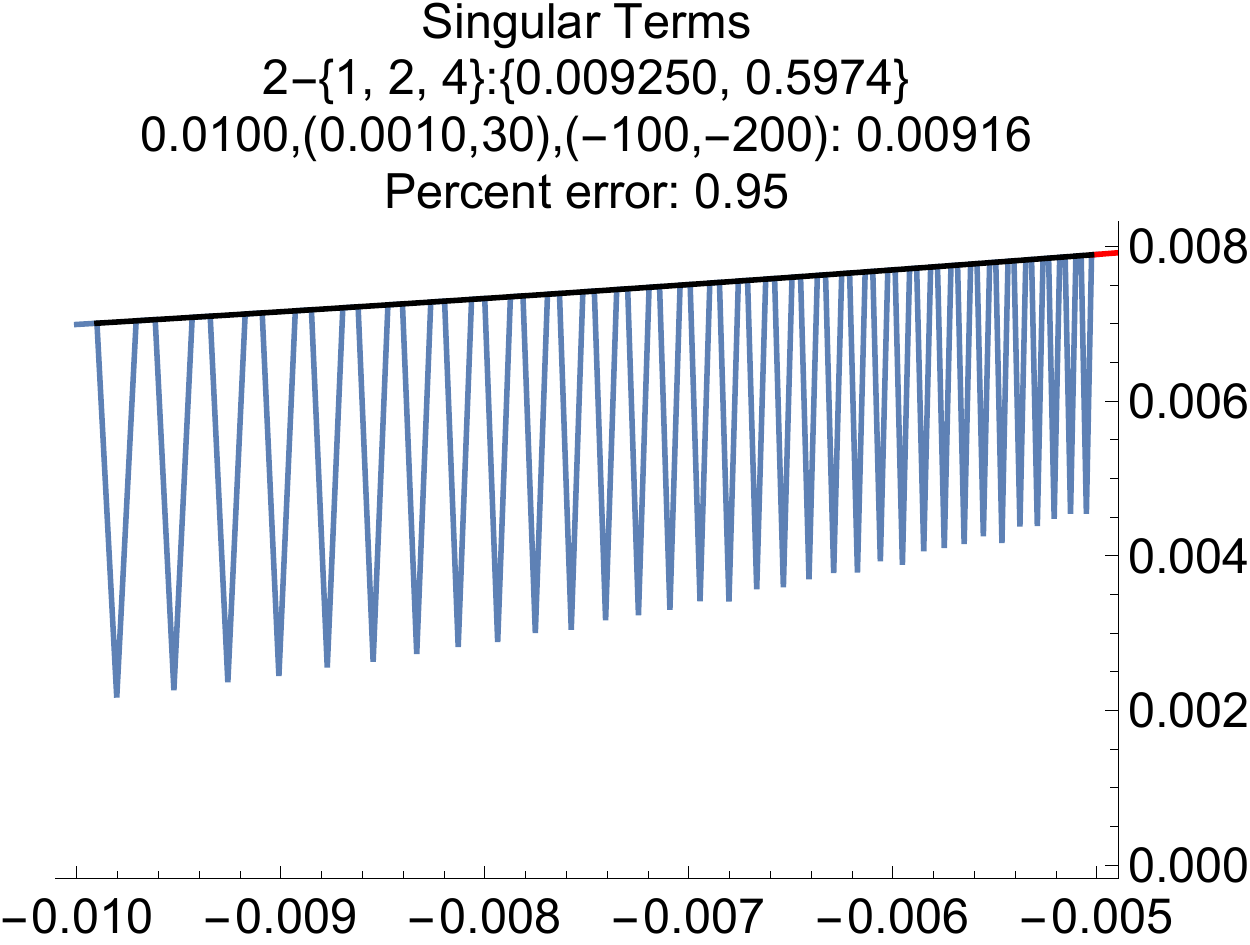}
			       }
     \hfill
     \subfloat[\label{ring2branches:2}]{%
       \includegraphics[width=0.4\textwidth]{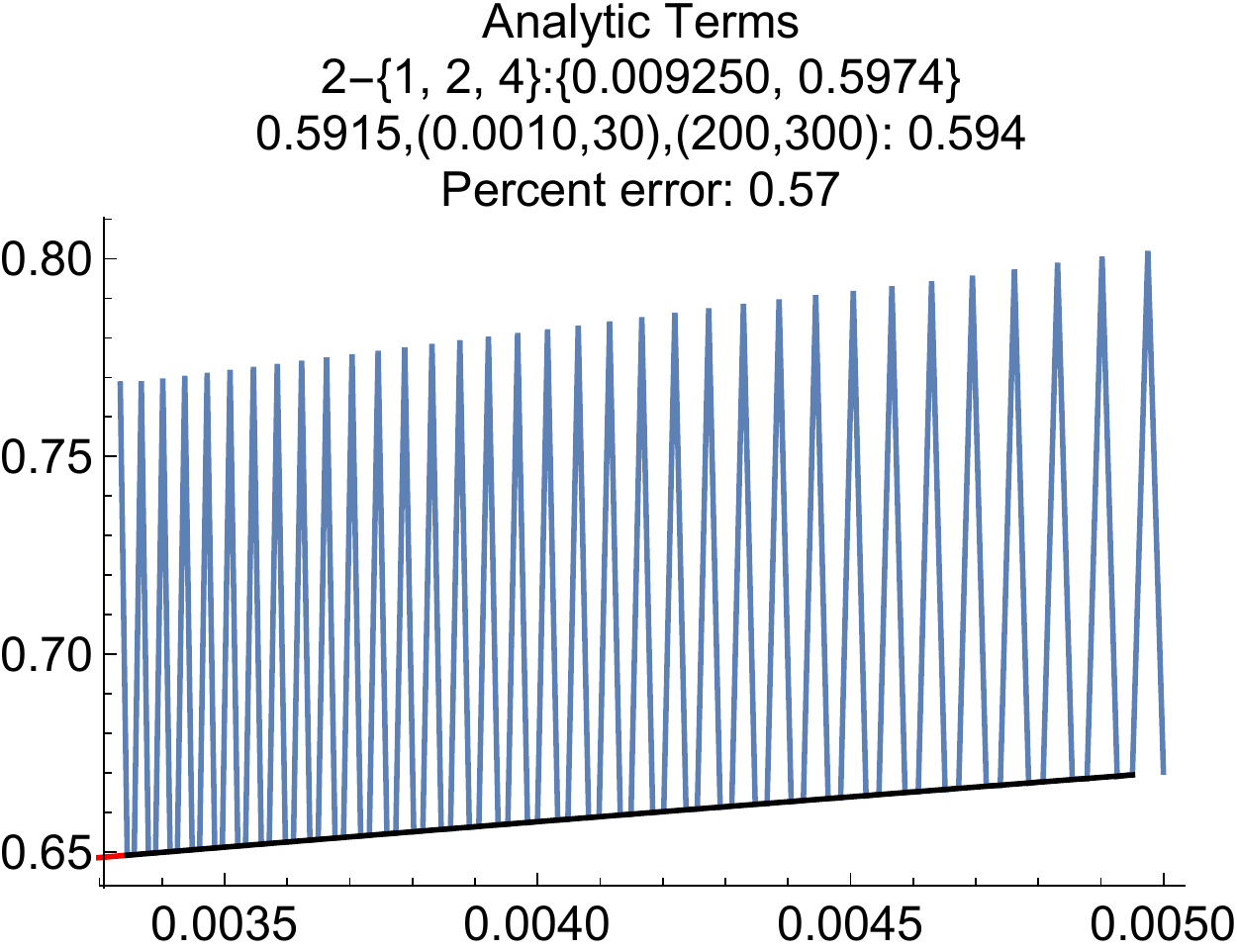}
     }
	     \hfill
      \caption{Root test results for annulus $2$, branch $(1,2,4)$}
     \label{figure:ring2branches}
   \end{figure}

 Figure \ref{figure:ring2branches} shows example plots of the Root Test results.  Note in the figure how the $\displaystyle\limsup$ of the Root Test data is used for extrapolating the singular data while $\displaystyle\liminf$ used for extrapolating the analytic terms.  The title for the singular plot (likewise for the analytic plot)gives the annulus and branch as $2-(1,2,4)$, the region of analyticity as determined by analytic continuation as $(0.00925,0.5974)$, radius of integration, step size and working precision of the calculations, as $0.01,(0.001,30)$, coefficients used as $(-200,-400)$, the extrapolated value returned by \texttt{FindFormula} as $0.009175$ as the percent error between the convergence value determined by analytic continuation and extrapolated Root Test.

\begin{table}[ht]
\caption{Root Test Results for Analytic Series} 
\centering 
\begin{tabular}{|c| c| c| c|c|c| c|c|} 
\hline
 \text{Annulus} &\text{Branch} & \text{Radius} & \text{Precision} & \text{Range} & \text{Actual } & \text{Est.}& \text{$\%$ Error} \\
\hline\hline 
 $1 $ & $ \{1,3\} $ & $ 0.009059 $ & $ 35 $ & $ \{100,200\} $ & $ 0.009150 $ & $ 0.009225 $ & $ 0.83 $ \\
 $1 $ & $ \{2\} $ & $ 0.6859 $ & $ 30 $ & $ \{50,100\} $ & $ 0.6929 $ & $ 0.6964 $ & $ 0.51 $ \\
 $1 $ & $ \{4\} $ & $ 0.009059 $ & $ 30 $ & $ \{50,100\} $ & $ 0.009150 $ & $ 0.009231 $ & $ 0.89 $ \\
 $2 $ & $ \{1,2,4\} $ & $ 0.5915 $ & $ 30 $ & $ \{200,300\} $ & $ 0.5974 $ & $ 0.5940 $ & $ 0.57 $ \\
 $3 $ & $ \{1,2\} $ & $ 0.6322 $ & $ 35 $ & $ \{100,200\} $ & $ 0.6325 $ & $ 0.6388 $ & $ 0.98 $ \\
 $3 $ & $ \{4\} $ & $ 0.8482 $ & $ 35 $ & $ \{100,200\} $ & $ 0.8507 $ & $ 0.8586 $ & $ 0.93 $ \\
 $4 $ & $ \{1\} $ & $ 0.8157 $ & $ 30 $ & $ \{50,100\} $ & $ 0.8175 $ & $ 0.8235 $ & $ 0.73 $ \\
 $4 $ & $ \{2\} $ & $ 0.6923 $ & $ 30 $ & $ \{50,100\} $ & $ 0.6929 $ & $ 0.6975 $ & $ 0.66 $ \\
 $5 $ & $ \{2,3\} $ & $ 0.8543 $ & $ 40 $ & $ \{300,600\} $ & $ 0.8559 $ & $ 0.8631 $ & $ 0.84 $ \\
 $6 $ & $ \{1\} $ & $ 0.8504 $ & $ 40 $ & $ \{200,400\} $ & $ 0.8507 $ & $ 0.8607 $ & $ 1.2 $ \\
 $7 $ & $ \{1\} $ & $ 0.8590 $ & $ 30 $ & $ \{200,300\} $ & $ 0.8591 $ & $ 0.8675 $ & $ 0.97 $ \\
 $7 $ & $ \{4\} $ & $ 0.8558 $ & $ 30 $ & $ \{50,100\} $ & $ 0.8559 $ & $ 0.8608 $ & $ 0.57 $ \\
 $8 $ & $ \{2,3\} $ & $ 0.8607 $ & $ 40 $ & $ \{200,300\} $ & $ 0.8607 $ & $ 0.8655 $ & $ 0.56 $ \\
 $8 $ & $ \{4\} $ & $ 0.8591 $ & $ 30 $ & $ \{50,100\} $ & $ 0.8591 $ & $ 0.8673 $ & $ 0.95 $ \\
 $9 $ & $ \{1,4\} $ & $ 0.8650 $ & $ 30 $ & $ \{200,300\} $ & $ 0.8727 $ & $ 0.8754 $ & $ 0.31 $ \\
 $10 $ & $ \{2,3\} $ & $ 0.8726 $ & $ 40 $ & $ \{300,600\} $ & $ 0.8727 $ & $ 0.8793 $ & $ 0.76 $ \\
 $11 $ & $ \{1\} $ & $ 0.9013 $ & $ 30 $ & $ \{50,100\} $ & $ 0.9016 $ & $ 0.9047 $ & $ 0.35 $ \\
 $11 $ & $ \{2,3,4\} $ & $ 0.9013 $ & $ 35 $ & $ \{100,200\} $ & $ 0.9016 $ & $ 0.8995 $ & $ 0.23 $ \\
 $12 $ & $ \{1,3,4,2\} $ & $ 0.9602 $ & $ 40 $ & $ \{300,600\} $ & $ 0.9608 $ & $ 0.9657 $ & $ 0.51 $ \\
 $13 $ & $ \{1,3,4,2\} $ & $ 0.9665 $ & $ 30 $ & $ \{200,300\} $ & $ 0.9666 $ & $ 0.9739 $ & $ 0.76 $ \\
 $14 $ & $ \{1,2\} $ & $ 1.190 $ & $ 40 $ & $ \{300,600\} $ & $ 1.192 $ & $ 1.200 $ & $ 0.65 $ \\
 $14 $ & $ \{3,4\} $ & $ 1.190 $ & $ 40 $ & $ \{300,600\} $ & $ 1.192 $ & $ 1.200 $ & $ 0.65 $ \\
 $15 $ & $ \{1,3\} $ & $ 1.400 $ & $ 40 $ & $ \{300,600\} $ & $ 1.403 $ & $ 1.405 $ & $ 0.15 $ \\
 $15 $ & $ \{2,4\} $ & $ 1.295 $ & $ 35 $ & $ \{100,200\} $ & $ 1.296 $ & $ 1.304 $ & $ 0.58 $ \\
 $16 $ & $ \{2\} $ & $ 1.303 $ & $ 30 $ & $ \{200,300\} $ & $ 1.303 $ & $ 1.311 $ & $ 0.60 $ \\
 $17 $ & $ \{2\} $ & $ 1.402 $ & $ 35 $ & $ \{100,200\} $ & $ 1.403 $ & $ 1.415 $ & $ 0.90 $ \\

\hline
\end{tabular}
\label{table:analRootTest} 
\end{table}

  It's important to note the plots filter out numerically imprecise values that are close to zero:  In the case of the singular data, we have the expression $\displaystyle 1/\left(\frac{c}{I(k,n)}\right)^{k/n}$.  Numerical precision will however result in small residual amounts that are actually zero which in turn results in very small terms in the plot which are neglected by $\displaystyle\limsup$.  A similar argument applies to the analytic terms.  However this does not mean the terms can be neglected but rather are simply not included in the calculations of the Root Test. 

Tables \ref{table:singRootTest} and \ref{table:analRootTest} summarize the results of the Root Tests.  Each table lists the branch, radius of integration, working precision of the numerical integrations, coefficients used in the Root Test, and the percent error between convergence data determined by analytic continuation and extrapolated data from the Root Test.   The singular table runs from annulus 2 through 18 since the first annulus has no lower limit.  And the analytic table runs from annulus 1 until it reaches a branch with no upper domain limit.  These are the branches with $\infty$ as the upper limit.  

\begin{figure}[!ht]
     \subfloat[\label{ring2brancheas:1}]{%
       \includegraphics[width=0.4\textwidth]{ring2branch124lowlim.pdf}
			       }
     \hfill
     \subfloat[\label{ring2branchesb:2}]{%
       \includegraphics[width=0.4\textwidth]{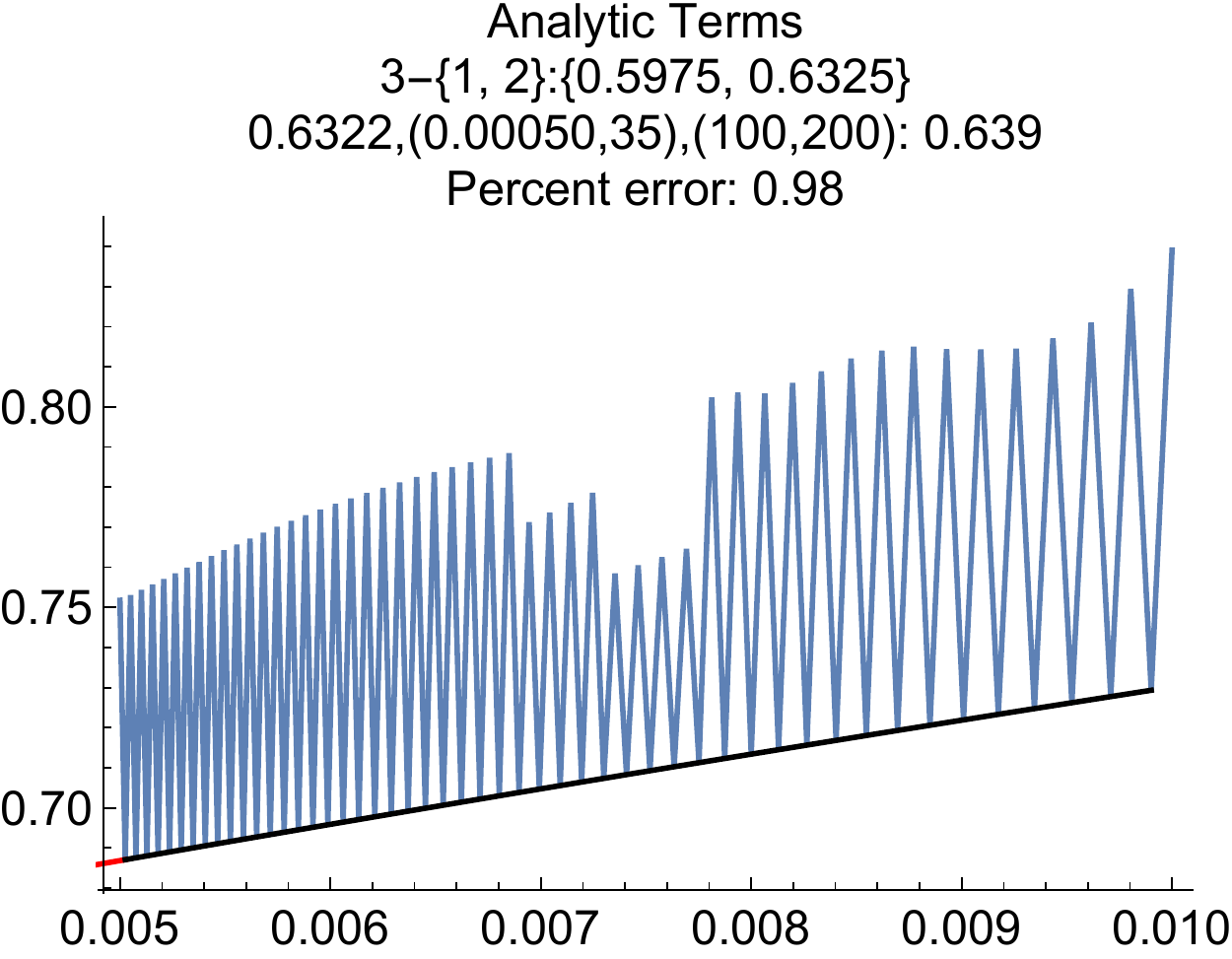}
     }
		\\
		 \subfloat[\label{ring2branchesc:1}]{%
       \includegraphics[width=0.4\textwidth]{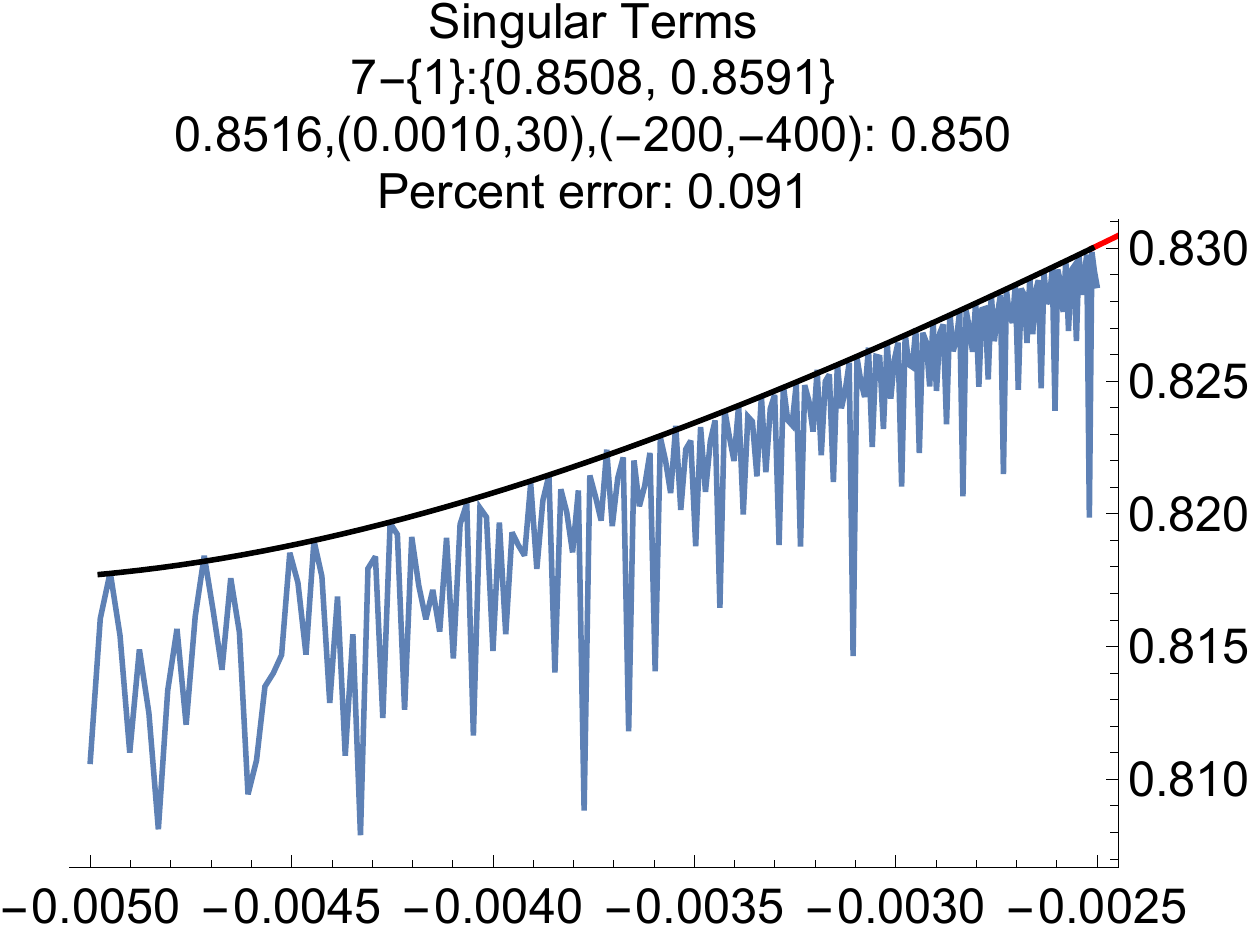}
			       }
     \hfill
     \subfloat[\label{ring2branchesd:2}]{%
       \includegraphics[width=0.4\textwidth]{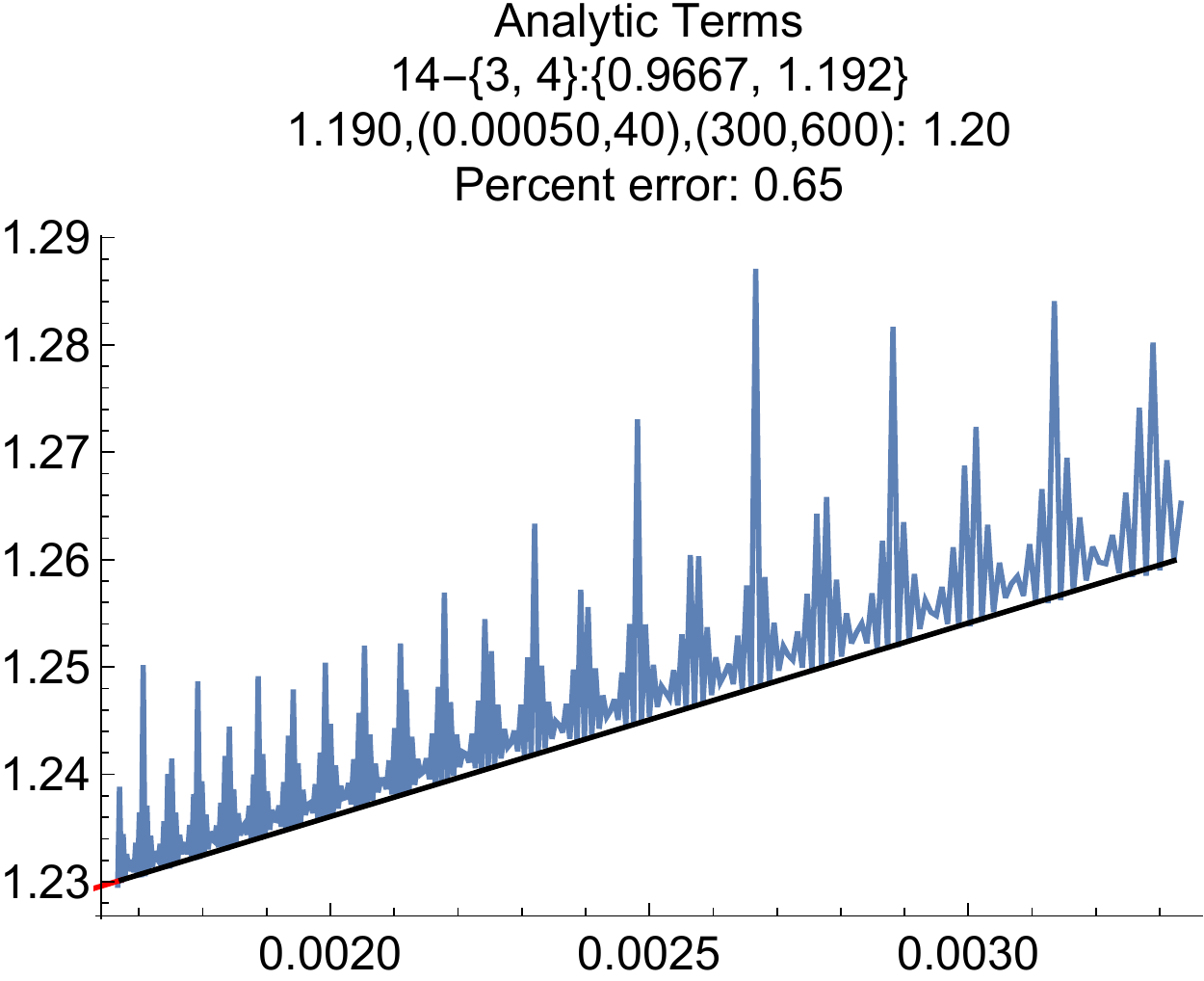}
     }		
	     \hfill
      \caption{Typical convergence morphologies for Root test}
     \label{figure:convergenceMorphologies}
   \end{figure}

Most results agreed to within $1\%$ error using the analysis parameters listed in each table.  Figure \ref{figure:convergenceMorphologies} shows typical morphologies of convergence data.

Several factors affected the results of the Root Test:
\begin{enumerate}
\item{Radius of Integration: }
The radius of integration affected the numerical precision of the results.  For example, in annulus 2 branch (1,2,4).  If $r=0.5$ then the numerical percision dropped quickly for the lower limit producing poor results.  However, if $r=0.001$, then the precision stabilized and the data converged nicely to the expected result.  Likewise, when the annular region extended beyond the unit circle, the term $r^k$ grows very rapidly so that we use a value close to the lower limit.
\item{A sufficient number of high-order coefficients were needed for the power series to settle down into a regular behavior: }
In some cases, only the first $100$ terms yielded a stable converging pattern.  In other cases, as many as $600$ terms were used.  And several test cases did not go below $1\%$ error even after analyzing $600$ terms at the level of precision used in this work.
\item{Results depended on the working precision and step size of the integration: }
 In this study, a maximum precision of $50$ digits and step size of $1/10000$ was used.
\item{Some coefficients of the power expansions are actually zero: }
Care must be taken to visually inspect the numerically-determined coefficients to decide if some terms are actually zero. Branch $(4)$ in annulus $16$ has an unbounded upper domain of convergence and all but the $a_0$ coefficient of the analytic terms are zero.  However, due to numerical precision of the integrations, we obtain residual non-zero values for the remaining terms which grow smaller with increasing working precision.
In the case of branch $(2,3)$ in annulus $18$, all of the analytic terms are zero.  This is because this $(2,3)$-cycle branch is singular at infinity.
\end{enumerate}
\bigskip
\section{Plotting the results}
 For readers interested in investigating equation \ref{equation:equation20} and it's branches further, the author has a website \cite{Milioto} with an interactive 3D viewer.  The viewer can be used to illustrate each branch.  See \href{http://jujusdiaries.com}{\color{blue}Algebraic Functions}.  
\section{Conclusions}
The branching geometry and domains determined by analytic continuation agreed well with results computed with the Root Test.  At a maximum of $50$ digits of precision and step size of $1/10000$, all results were in agreement below $2\%$ error.  These results suggests the analytic continuation method described above can be used to successfully determine the annular branching domains of algebraic functions of low degree, and also demonstrated the Root Test can successfully approximate the convergence domains of the associated power expansions.  Care should be taken to inspect the integration results manually to confirm the integrity of the data especially since the integrands of Equation \ref{eqn5} become highly oscillatory with increasing $k$ and therefore increasingly difficult to compute accurately. 


\begin{thebibliography}{10}  

\bibitem {Bliss} Bliss, Gilbert A. \textit{Algebraic Functions}. New York: Dover Publications, Inc., 2004.

\bibitem {Brown} Brown, James and Ruel Churchill. \textit{Complex Variables and Applications}. New York:  McGraw Hill, 2004

\bibitem {Chud} Chudnovsky, D.V. and G.V. Chudnovsky. ``On Expansion of Algebraic Functions in Power and Puiseux Series''. Journal of Complexity \textbf{2}, 271-294 (1986).

\bibitem {Kung} Kung, H.T. and J. Traub, ``All Algebraic Functions can be Computed Fast''. J. Assoc. Comput. Mach. \textbf{25}, 245-260.

\bibitem {Marku} Markushevich,A.I.,1967.\textit{Theory of Functions of a Complex Variable.Vol.III.} PrenticeHall,
Englewood Cli?s, N. J.

\bibitem {Marsden} Marsden, Jerrold and Michael Hoffman. \textit{Basic Complex Analysis}. New York: W.H Freeman and Company, 1999.
\bibitem {Milioto} Milioto,Dominic C. (2018, Dec. 8). Algebraic Functions, Retrieved from http://jujusdiaries.com.
\bibitem {Nowak} Nowak, Krzysztof.\ \textit{Some Elementary Proofs of Puiseux's Theorem}.  Universitatis Iagellonicae ACTA Mathematica, Fasciculus XXXVIII, 2000. 
\bibitem {Walker} Walker, Robert J.  \textit{Algebraic Curves}. Princeton:  Princeton University Press, 1956.

\bibitem {Willis} Willis, Nicholas J., Didier, Annie K., Sonnanburg, Kevin M. \textit{How to Compute a Puiseux Expansion},
arXiv: 0807.4674.1 [math.AG] 29 July, 2008



\end{thebibliography}
\end{document}